\documentclass[a4paper,11pt]{article}
\usepackage{amsmath}
\usepackage{amssymb}
\usepackage{enumerate}
\usepackage{theorem}
\usepackage{array}
\usepackage{bm}
\newtheorem{theorem}{Theorem}[section]
\newtheorem{lemma}[theorem]{Lemma}
\newtheorem{corollary}[theorem]{Corollary}
\newtheorem{proposition}[theorem]{Proposition}
\theorembodyfont{\rmfamily}
\newtheorem{definition}[theorem]{Definition}
\newtheorem{notation}[theorem]{Notation}
\newtheorem{remark}[theorem]{Remark}

\newtheorem{annotation}{Annotation}

\newtheorem{conjecture}{Conjecture}

\newtheorem{setup}{Set-up}
\newcommand{\proof}{\noindent \mbox{\em Proof.\hspace*{2mm}}}
\newcommand{\qed}{\hfill \mbox{$  \Box $}}
\DeclareMathOperator{\Sing}{Sing}

\newcommand{\gyokan}{\vskip 6pt}
\title{Number of points of a nonsingular hypersurface in
an odd-dimensional projective space}
\author{
Masaaki Homma
\thanks{Partially supported by JSPS KAKENHI Grant Number JP15K04829.}
\\
 Department of Mathematics and Physics\\
Kanagawa University\\
Hiratsuka 259-1293, Japan\\
homma@kanagawa-u.ac.jp
\and
Seon Jeong Kim
\thanks{Partially supported by Basic Science Research Program through the National Research Foundation of Korea(NRF) 
funded by the Ministry of Education (2016R1D1A1B01011730).
}\\
 Department of Mathematics and RINS\\
Gyeongsang National University\\
Jinju 660-701, Korea \\
skim@gnu.kr
}
\date{}

\begin{document}
\maketitle
\begin{abstract}
The numbers of $\mathbb{F}_q$-points of nonsingular hypersurfaces
of a fixed degree in an odd-dimensional projective space are investigated,
and an upper bound for them is given.
Also we give the complete list of nonsingular hypersurfaces
each of which realizes the upper bound.
This is a natural generalization of our previous study of surfaces in
projective $3$-space.
\end{abstract}

\section{Introduction}
Several years ago, we established the elementary bound for the numbers of
$\mathbb{F}_q$-points of hypersurfaces of projective $n$-space $\mathbb{P}^n$
with $n \geq 3$ \cite{hom-kim2013b},
and later gave the complete list of surfaces in $\mathbb{P}^3$
whose number of $\mathbb{F}_q$-points reached this bound 
\cite{hom-kim2015, hom-kim2015online}.
Recently Tironi extended this list for hypersurfaces in $\mathbb{P}^n$
\cite{tir2014preprint}.
Although surfaces appeared in the list are nonsingular,
hypersurfaces appeared in the extended list with $n>3$ are
cones over those surfaces except when the degree of the hypersurface is
$q+1$.
Therefore if we restrict our investigation within nonsingular hypersurfaces,
we can expect a tighter bound than the elementary bound.

\gyokan

We fix a finite field $\mathbb{F}_q$ of $q$ elements.
The number of $\mathbb{F}_q$-points of the projective $m$-space
is denoted by $\theta_q(m)$, that is,
$\theta_q(m)=\sum_{\nu=0}^{m} q^{\nu}$.
A closed subscheme $\mathcal{X}$ in $\mathbb{P}^m$
over the algebraic closure of $\mathbb{F}_q$
is said to be defined over $\mathbb{F}_q$ if
the homogeneous ideal of $\mathcal{X}$ is generated by
polynomials $f_1(X_0, \dots , X_m), \dots , f_s(X_0, \dots , X_m)$
in $\mathbb{F}_q[X_0, \dots , X_m]$.
An $\mathbb{F}_q$-point $(a_0, \dots , a_m)$ of $\mathbb{P}^m$
is said to be an $\mathbb{F}_q$-point of $\mathcal{X}$
if $f_1(a_0, \dots , a_m)= \dots = f_s(a_0, \dots , a_m)=0$,
namely we do not care the point is a multiple point or not in
$\mathcal{X}$.
The set of $\mathbb{F}_q$-points of $\mathcal{X}$
is denoted by $\mathcal{X}(\mathbb{F}_q)$
and the cardinality of this set by $|\mathcal{X}(\mathbb{F}_q)|$
or $N_q(\mathcal{X})$.
We frequently use the notation
$\{ f_1 = \dots =f_m =0\}$ for the scheme $\mathcal{X}$.

Geometric structure of $\mathcal{X}$, for example, nonsingularity
or irreducibility,
is normally (and also in this article) considered over
the algebraic closure $\overline{\mathbb{F}}_q$
of $\mathbb{F}_q$,
but we are just interesting in the set-theoretical counting of
$\mathcal{X}(\mathbb{F}_q)$.

The purpose of this article is to show the following theorem.

\begin{theorem}\label{maintheorem}
Let $n$ be an odd integer at least $3$.
If $X$ is a nonsingular hypersurface of degree $d\geq 2$
in $\mathbb{P}^n$ defined over $\mathbb{F}_q$.
Then
\[
N_q(X) \leq \theta_q \! \left(\frac{n-1}{2}\right) \cdot
\left((d-1)q^{\frac{n-1}{2}} + 1\right),
\]
and equality holds if and only if either
\begin{enumerate}[{\rm (i)}]
\item $d=2$ and $X$ is the nonsingular hyperbolic quadric hypersurface,
that is, $X$ is projectively equivalent over $\mathbb{F}_q$ to the hypersurface
\[
\sum_{i=0}^{\frac{n-1}{2}} X_{2i}X_{2i+1} =0 \, \rm{ ; } \text{ or}
\]
\item $d=\sqrt{q} +1$ where $q$ is square, and $X$ is
a nonsingular Hermitian hypersurface, that is,
$X$ is projectively equivalent over $\mathbb{F}_q$ to the hypersurface
\[
\sum_{i=0}^{\frac{n-1}{2}} 
  \left( X_{2i}^{\sqrt{q}}X_{2i+1}  + X_{2i}X_{2i+1}^{\sqrt{q}}\right)= 0 
 \,  \rm{;} \text{ or}
\]
\item $d=q+1$ and $X$ is a nonsingular $\mathbb{P}^n$-filling hypersurface
over $\mathbb{F}_q$, that is,
$X$ is projectively equivalent over $\mathbb{F}_q$ to the hypersurface
\[
\sum_{i=0}^{\frac{n-1}{2}} 
   \left( X_{2i}^q X_{2i+1}  - X_{2i}X_{2i+1}^q\right) =0 
   \text{.}
\]
\end{enumerate}
\end{theorem}

We prove this by induction on $n$, so $n=3$ is the first step of the induction,
which was already showed in \cite[Theorem~1]{hom-kim2015}:

\begin{theorem}\label{nequal3}
Let $X$ be a surface of degree $d$ in $\mathbb{P}^3$ over $\mathbb{F}_q$
without $\mathbb{F}_q$-plane components.
Then
$
N_q(X) \leq \theta_q(1) \cdot ((d-1)q +1),
$
and 
equality holds if and only if 
the degree $d$ is either
$2$ or $\sqrt{q} +1$ {\rm(}when $q$ is a square{\rm)} or $q+1$
and
the surface $X$ is projectively equivalent to
one of the following surfaces over ${\Bbb F}_q$
according to the degree{\rm :}
\begin{enumerate}[{\rm (i)}]
\item
$X_0X_1 + X_2X_3=0$  if $d=2${\rm ;}
\item
$X_0^{\sqrt{q}}X_1+ X_0X_1^{\sqrt{q}} +X_2^{\sqrt{q}}X_3+X_2X_3^{\sqrt{q}}=0$
if $d= \sqrt{q}+1${\rm ;}
\item $X_0^{q}X_1 - X_0X_1^{q} + X_2^{q}X_3 - X_2X_3^{q} =0$
if $d=q+1$.
\end{enumerate}
\end{theorem}
\begin{remark}
\begin{enumerate}[(i)]
\item The assumption that
$X$ has no ${\Bbb F}_q$-plane components
in the above theorem is milder than the nonsingularity of $X$
if $\deg X \geq 2$.
\item Equations in the above theorem and those
in \cite[Theorem~1]{hom-kim2015}
are seemingly different.
But one can easily confirm that in each degree those equations are projectively equivalent
over ${\Bbb F}_q$ to each other.
\end{enumerate}
\end{remark}
\section{Preliminary}
This section is a mixture of facts that are mostly independent of one another,
but necessary to our proof.

We keep roman letters $X, Y, Z$ for particular varieties for later use.
In this section, varieties or schemes are denoted by calligraphic letters
$\mathcal{X}, \mathcal{Y}, \mathcal{Z}$ etc. 

\subsection{A necessary condition of a hypersurface to be nonsingular}
\begin{lemma}\label{linearinnonsingular}
Let $\mathcal{X}$ be a hypersurface of degree $\geq 2$ in $\mathbb{P}^m$
over an algebraically closed field, and $\mathcal{L}$ a linear subspace
of $\mathbb{P}^m$
which is contained in $\mathcal{X}$.
If $\mathcal{X}$ is nonsingular,
then $\dim \mathcal{L} \leq \lfloor\frac{m-1}{2} \rfloor$.
Here the symbol $\lfloor \frac{m-1}{2} \rfloor$ denotes the integer part
of $\frac{m-1}{2}$.
\end{lemma}
\proof
Let $r = \dim \mathcal{L}$.
Choose the coordinates $X_0, \dots , X_m$ of $\mathbb{P}^m$
so that $\mathcal{L}$ is defined by
$X_0 = X_1 = \dots = X_{m-r-1}=0.$
Since $\mathcal{L} \subset \mathcal{X}$,
the equation of $\mathcal{X}$ is of the form
\[
F(X_0, \dots , X_m) = \sum_{i=0}^{m-r-1} f_i(X_0, \dots , X_m) X_i=0.
\]
Note that each homogeneous polynomial $f_i(X_0, \dots , X_m)$
is not constant because $\deg \mathcal{X} \geq 2$.
Consider the simultaneous equations
\[
F = \cfrac[l]{\partial F}{\partial X_0}= \dots 
 = \cfrac[l]{\partial F}{\partial X_{m}} =0,
\]
more explicitly:
\begin{equation}
 \left\{
 \begin{split}\label{partialderivative}
  &F = \sum_{i=0}^{m-r-1} f_i X_i=0\\
  &\cfrac[l]{\partial F}{\partial X_0} 
         = \sum_{i=0}^{m-r-1} \cfrac[l]{\partial f_i}{\partial X_0}X_i 
                   +f_0  =0\\
           &\qquad \vdots \\
  &\cfrac[l]{\partial F}{\partial X_{m-r-1}} 
      = \sum_{i=0}^{m-r-1} \cfrac[l]{\partial f_i}{\partial X_{m-r-1}}X_i 
             +f_{m-r-1} =0\\  
     & \cfrac[l]{\partial F}{\partial X_{m-r}} 
      = \sum_{i=0}^{m-r-1} \cfrac[l]{\partial f_i}{\partial X_{m-r}}X_i =0\\
            &\qquad \vdots \\
  &\cfrac[l]{\partial F}{\partial X_{m}} 
      = \sum_{i=0}^{m-r-1} \cfrac[l]{\partial f_i}{\partial X_{m}}X_i=0 .
 \end{split}
 \right.
\end{equation}
We may view $\{ X_{m-r}, \dots , X_m \}$
as a system of coordinates of $\mathcal{L} = \mathbb{P}^r$.
Suppose $m-r \leq r.$
Then the simultaneous $m-r$ equations
\begin{equation}\label{simultaneous}
 \left\{
 \begin{split}
  &f_0(0, \dots , 0, X_{m-r}, \dots , X_m ) = 0 \\
                &\qquad \vdots \\
  &f_{m-r-1}(0, \dots ,0,  X_{m-r}, \dots , X_m ) = 0 
 \end{split}
 \right.
\end{equation}
has a solution $(\alpha_{m-r}, \dots , \alpha_m)$
in $\mathbb{P}^r$.
Hence the point $(0, \dots , 0, \alpha_{m-r}, \dots , \alpha_m)$
in $\mathcal{L} \subset \mathcal{X}$
is a solution of (\ref{partialderivative}),
which must be a singular point of $\mathcal{X}$.
Therefore we have $m-r > r$
if $\mathcal{X}$ is nonsingular.
\qed

\subsection{Segre-Serre-S{\o}rensen bound}
Without any restrictions on a hypersurface over $\mathbb{F}_q$,
the best bound was obtained by Serre \cite{ser1991},
which is a generalization of Segre's old result
for plane curves \cite{seg1959}.
S{\o}rensen \cite{sor1994} also proved the same inequality as Serre's.
\begin{lemma}[Segre-Serre-S{\o}rensen]\label{theoremSSS}
Let $\mathcal{X} \subset \mathbb{P}^m$
be a hypersurface of degree $d$
defined over $\mathbb{F}_q$.
Then
$N_q(\mathcal{X}) \leq dq^{m-1} + \theta_q(m-2)$.
Moreover, when ``$m=2$" or ``$m>2$ and $d \leq q$",
equality holds if and only if
there are $d$ hyperplanes 
$\mathcal{L}_1, \dots , \mathcal{L}_d$
over $\mathbb{F}_q$
such that
$\mathcal{X} = \cup_{i=1}^{d} \, \mathcal{L}_i$
and $\cap_{i=1}^{d} \, \mathcal{L}_i$
is of dimension $m-2$.
\end{lemma}
\proof
See \cite[II \S 6 Observation IV]{seg1959} for ``$m=2$",
\cite{ser1991} for ``$m>2$".

\begin{notation}
For a variety $\mathcal{X}$, $\Sing \mathcal{X}$
denotes the locus of singular points.
\end{notation}
In Lemma~\ref{theoremSSS}, 
$\Sing \mathcal{X} = \cap_{i=1}^{d} \, \mathcal{L}_i.$
Actually, the following lemma holds.
\begin{lemma}\label{singularlocus}
Let $\mathcal{X}$ be a hypersurface in $\mathbb{P}^m$.
If $\mathcal{X}$ splits into hyperplanes {\rm :}
$\mathcal{X} = \cup_{i=1}^{d} \mathcal{L}_i,$
then $\Sing \mathcal{X} = \cup_{i < j} (\mathcal{L}_i \cap \mathcal{L}_j).$
\end{lemma}
\proof
Let $g_i = \sum_{j=0}^m a_{ij}X_j=0$
be the linear equation of $\mathcal{L}_i$.
So $\mathcal{X}$ is defined by $G=\prod_{i=1}^{d} g_i=0$.
Then
\[
\cfrac[l]{\partial G}{\partial X_{\nu}}
= \sum_{i=1}^{d} a_{i\nu} \prod_{\begin{subarray}{c}
                       l \text{ with}\\
                       l \neq i
                      \end{subarray}} g_l.
                      \]
If $(u_0, \dots , u_m) \in \mathcal{L}_{\alpha}
          \cap \mathcal{L}_{\beta}$,
then 
$\cfrac[l]{\partial G}{\partial X_{\nu}}(u_0, \dots , u_m)=0$
because $g_{\alpha}$ or $g_{\beta}$ appears in each term of
$\cfrac[l]{\partial G}{\partial X_{\nu}}$.
Hence $\mathcal{L}_{\alpha}
          \cap \mathcal{L}_{\beta} \subset \Sing \mathcal{X}.$
Conversely if 
$(u_0, \dots , u_m) \in \mathcal{X} \setminus 
                    \cup_{i<j} \left( \mathcal{L}_i
          \cap \mathcal{L}_j\right) ,$
there is a unique hyperplane $\mathcal{L}_{\nu}$
which contains the point $(u_0, \dots , u_m)$.
Hence $\cfrac[l]{\partial G}{\partial X_{\nu}}(u_0, \dots , u_m)
=a_{\alpha \nu} \prod_{l \neq \alpha} g_l(u_0, \dots , u_m)$,
is nonzero for some $\nu$.
Hence $\mathcal{X} \setminus 
                    \cup_{i<j} \left( \mathcal{L}_i
          \cap \mathcal{L}_j \right) \subset \mathcal{X} \setminus \Sing \mathcal{X}$.
\qed

\gyokan

We frequently use the latter half of Segre-Serre-S{\o}rensen's lemma
(\ref{theoremSSS}).
For the convenience of readers,
we reformulate the necessary part with a small generalization
and give its proof.

\begin{lemma}\label{dL}
Let $\mathcal{L}_1, \dots , \mathcal{L}_d$ be distinct
linear subspaces over $\mathbb{F}_q$
in $\mathbb{P}^m$ such that
\begin{enumerate}[{\rm (i)}]
  \item $\dim\, \mathcal{L}_1 = \dots = \dim\, \mathcal{L}_d = k$, and
  \item $\dim \, \bigcap_{i=1}^d \mathcal{L}_i= k-1$.
\end{enumerate}
Then $N_q(\mathcal{L}_1 \cup \dots \cup \mathcal{L}_d)
       = d q^k + \theta_q(k-1).$
\end{lemma}
\proof
Let $\Lambda = \bigcap_{i=1}^{d} \mathcal{L}_i$.
From the assumptions, $\mathcal{L}_i \cap \mathcal{L}_j= \Lambda$
if $i \neq j$.
Therefore
\[
(\mathcal{L}_1 \cup \dots \cup \mathcal{L}_d)(\mathbb{F}_q)
 = \Bigl( \coprod_{i=1}^d ( \mathcal{L}_i(\mathbb{F}_q) \setminus 
                                      \Lambda(\mathbb{F}_q))\Bigr)
                             \coprod \Lambda(\mathbb{F}_q),
\]
where the symbol $\coprod$ means taking the disjoint union.
Hence
$
N_q(\mathcal{L}_1 \cup \dots \cup \mathcal{L}_d)
= d (\theta_q(k) - \theta_q(k-1)) + \theta_q(k-1).
$
\qed

\gyokan

The next lemma is also useful.

\begin{lemma}\label{nonsingularity_ancestor}
Let $\mathcal{X}$ be a hypersurface of $\mathbb{P}^m$,
and $\mathcal{S}$ a linear subspace of $\mathbb{P}^m$
such that $\mathcal{S} \not\subset \mathcal{X}.$
If a point $Q \in \mathcal{S} \cap \mathcal{X}$
is nonsingular in $\mathcal{S} \cap \mathcal{X}$,
then $Q$ is also nonsingular in $\mathcal{X}$.
\end{lemma}
\proof
We assume that $\mathcal{S}=\{X_0 = \dots = X_s =0\}$
and $Q=(0, \dots , 0,1)$.
Use affine coordinates
$x_0 = \frac{X_0}{X_m}, \dots , x_{m-1} = \frac{X_{m-1}}{X_m}$.
Let $f(x_0, \dots , x_{m-1}) = f_1 + f_2 + \dots + f_d=0$
be the local equation of $\mathcal{X}$ around $Q$,
where $f_i=f(x_0, \dots , f_{m-1})$ is the homogeneous part of degree
$i$ of $f$.
Since $Q$ is nonsingular in $\mathcal{S} \cap \mathcal{X}$,
$f_1(0, \dots , 0, x_{s+1}, \dots , x_{m-1})$ is nontrivial.
Hence so is $f_1(x_1, \dots , x_{m-1})$.
\qed

\subsection{Cone lemma}
\begin{lemma}\label{triviallemma}
Let $f(X_0, \dots , X_m)$ be a homogeneous polynomial
over $\mathbb{F}_q$ of degree $d \leq q$.
If $f(a_0, \dots , a_m)=0$ for any
$(a_0, \dots , a_m) \in \mathbb{F}_q^n,$
then $f$ is the zero polynomial.
\end{lemma}
\proof
Suppose $f$ is nontrivial,
then it defines a hypersurface $\mathcal{X}$ of degree $d$ in $\mathbb{P}^m$.
By the lemma of Segre-Serre-S{\o}rensen (\ref{theoremSSS}),
\begin{align*}
 N_q(\mathcal{X}) & \leq dq^{m-1} + \theta_q(m-2)  \\
                  & \leq q^m +  \theta_q(m-2) 
                    <  \theta_q(m) \  \, \text{ if  \ $d \leq q$,} 
\end{align*}
however, $\mathcal{X}(\mathbb{F}_q) = \mathbb{P}^m(\mathbb{F}_q)$
by the assumption; these contradict each other.
Hence the polynomial must be trivial.
\qed

The following fact, which will be referred to as ``cone lemma",
is a bridge between point-counting and geometry.
\begin{proposition}\label{conelemma}
Let $\mathcal{X}$ be a hypersurface over $\mathbb{F}_q$ of degree $d$
with $2 \leq d \leq q$
in $\mathbb{P}^m$, and $\mathcal{L}=\mathbb{P}^{m-k}$
an $\mathbb{F}_q$-linear subvariety of $\mathbb{P}^m$ of dimension $m-k$,
where $3 \leq k \leq m$.
Let $\mathcal{M}=\mathbb{P}^{k-1}$ be another $\mathbb{F}_q$-linear subvariety
of $\mathbb{P}^m$ of dimension $k-1$
such that $\mathcal{L} \cap \mathcal{M} = \emptyset$,
and $\mathcal{Y}$ a hypersurface of degree $d$ in
$\mathbb{P}^{k-1} = \mathcal{M}$ over $\mathbb{F}_q$.
Suppose that
\begin{gather*}
  N_q(\mathcal{Y}) > (d-1)q^{k-2} + \theta_q(k-3), \\
  \mathcal{X} \supset \mathcal{Y}, \text{ and}\\
  \mathcal{X} \supset (\mathbb{P}^{m-k} \ast \mathcal{Y})(\mathbb{F}_q),
\end{gather*}
where $\mathbb{P}^{m-k} \ast \mathcal{Y}$ denotes the cone over 
$\mathcal{Y}$ with center $\mathbb{P}^{m-k}= \mathcal{L}.$
Then $\mathcal{X} = \mathbb{P}^{m-k} \ast \mathcal{Y}$.
\end{proposition}
\proof
Choose coordinates
$
X_0, \dots , X_{k-1}, X_k, \dots , X_m
$
of $\mathbb{P}^m$ so that
$\mathcal{L}=\mathbb{P}^{m-k}$ is defined by
$X_0 = \dots = X_{k-1}=0$, and
$\mathcal{M}=\mathbb{P}^{k-1}$
by $X_k = \dots = X_m =0.$
Let
\[
F(X_0, \dots , X_m) =
\sum_{
\begin{subarray}{c}
\bm{e}=(e_0, \dots , e_m) \\
\text{with}\\
e_0 + \dots + e_m =d
\end{subarray}
}
\alpha_{\bm{e}}X_0^{e_0} \cdots X_m^{e_m}
=0
\]
be the equation of $\mathcal{X}$.
Note that the polynomial $F$
can be rewritten as
\begin{multline}\label{eqcalX}
F(X_0, \dots , X_m) =\\
\sum_{\mu=0}^{d} \sum_{
\begin{subarray}{c}
(e_k, \dots , e_m) \\
\text{with}\\
\sum_k^m e_j=\mu
\end{subarray}
}
\Bigl(
\sum_{
\begin{subarray}{c}
(e_0, \dots , e_{k-1}) \\
\text{with}\\
\sum_0^{k-1} e_i =d-\mu
\end{subarray}
}
\alpha_{(e_0, \dots , e_{k-1}, e_k, \dots , e_m)}
X_0^{e_0}\cdots X_{k-1}^{e_{k-1}}
\Bigr) X_k^{e_k}\cdots X_m^{e_m}.
\end{multline}
Let $(0, \dots , 0, b_k, \dots b_m) \in \mathcal{L}(\mathbb{F}_q)$
and $(a_0, \dots , a_{k-1}, 0, \dots, 0) \in \mathcal{Y}(\mathbb{F}_q)$.
Since 
$\mathcal{X} \supset (\mathbb{P}^{m-k} \ast \mathcal{Y})(\mathbb{F}_q),$
\begin{equation}\label{ruling}
(ta_0, \dots , ta_{k-1},sb_k, \dots sb_m) \in \mathcal{X}(\mathbb{F}_q)
\end{equation}
for any $(s,t) \in \mathbb{P}^1(\mathbb{F}_q).$
Substitute (\ref{ruling}) for $F(X_0, \dots , X_m)$, then by (\ref{eqcalX})
\begin{equation}\label{eqsandt}
\sum_{\mu=0}^{d} t^{d-\mu}s^{\mu}\sum_{
\begin{subarray}{c}
(e_k, \dots , e_m) \\
\text{with}\\
\sum_k^m e_j=\mu
\end{subarray}
}
\Bigl(
\sum_{
\begin{subarray}{c}
(e_0, \dots , e_{k-1}) \\
\text{with}\\
\sum_0^{k-1} e_i =d-\mu
\end{subarray}
}
\alpha_{(e_0, \dots , e_{k-1}, e_k, \dots , e_m)}
a_0^{e_0}\cdots a_{k-1}^{e_{k-1}}
\Bigr) b_k^{e_k}\cdots b_m^{e_m} = 0
\end{equation}
for any $(s,t) \in \mathbb{P}^1(\mathbb{F}_q).$
Since $d \leq q $ but $|\mathbb{P}^1(\mathbb{F}_q)|=q +1$, 
all coefficients of the polynomial (\ref{eqsandt})
in $s$ and $t$ are $0$.
Hence
\begin{equation}\label{eqandb}
\sum_{
\begin{subarray}{c}
(e_k, \dots , e_m) \\
\text{with}\\
\sum_k^m e_j=\mu
\end{subarray}
}
\Bigl(
\sum_{
\begin{subarray}{c}
(e_0, \dots , e_{k-1}) \\
\text{with}\\
\sum_0^{k-1} e_i =d-\mu
\end{subarray}
}
\alpha_{(e_0, \dots , e_{k-1}, e_k, \dots , e_m)}
a_0^{e_0}\cdots a_{k-1}^{e_{k-1}}
\Bigr) b_k^{e_k}\cdots b_m^{e_m} = 0
\end{equation}
for any
$(a_0, \dots , a_{k-1}, 0, \dots, 0) \in \mathcal{Y}(\mathbb{F}_q)$
and
$(b_k, \dots , b_m) \in \mathbb{P}^{m-k}(\mathbb{F}_q).$

First we fix the element $(a_0, \dots , a_{k-1}) \in \mathcal{Y}(\mathbb{F}_q)$
and view (\ref{eqandb}) as a polynomial
with variables $(b_k, \dots , b_m)$.
Since the degree of the polynomial (\ref{eqandb}) in $(b_k, \dots , b_m)$
is $\mu$ ($\leq d \leq q$)
and it is $0$ for any $(b_k, \dots , b_m) \in \mathbb{P}^{m-k}(\mathbb{F}_q),$
it must be a zero polynomial by (\ref{triviallemma}),
that is
\[
\sum_{
\begin{subarray}{c}
\bm{e}=(e_0, \dots , e_{k-1}) \\
\text{with}\\
e_0 + \dots + e_{k-1} =d-\mu
\end{subarray}
}
\alpha_{(\bm{e}, e_k, \dots , e_m)}
a_0^{e_0}\cdots a_{k-1}^{e_{k-1}}
= 0
\]
for any $(e_k, \dots , e_m)$
with $e_k + \dots + e_m = \mu$.
Hence for each $\bm{e}'=(e_k, \dots , e_m)$
the hypersurface $\mathcal{Y}_{\bm{e}'}$
defined by
\[
\sum_{
\begin{subarray}{c}
\bm{e}=(e_0, \dots , e_{k-1}) \\
\text{with}\\
e_0 + \dots + e_{k-1} =d-\mu
\end{subarray}
}
\alpha_{(\bm{e}, \bm{e}')}
X_0^{e_0}\cdots X_{k-1}^{e_{k-1}}
= 0
\]
in $\mathbb{P}^{k-1}= \mathcal{M}$
contains $\mathcal{Y}(\mathbb{F}_q)$.
If the polynomial is nontrivial (of degree $d-\mu$), then
$N_q(\mathcal{Y}_{\bm{e}'}) \leq (d-\mu)q^{k-2}+\theta_q(k-3)$
by Lemma~\ref{theoremSSS}.
On the other hand,
$
N_q(\mathcal{Y}) > (d-1)q^{k-2} + \theta_q(k-3)
$
by the assumption.
Hence if $\mu \geq 1$,
this polynomial must be trivial.
Therefore $\mathcal{X}$
is a cone of the hypersurface
\[
\sum_{
\begin{subarray}{c}
(e_0, \dots , e_{k-1})\\
\text{with}\\
e_0 + \dots + e_{k-1} = d
\end{subarray}
}
\alpha_{(e_0, \dots , e_{k-1}, 0, \dots , 0)}X_0^{e_0}\cdots X_{k-1}^{e_{k-1}}
=0
\]
of $\mathbb{P}^{k-1}$, which is the equation of $\mathcal{Y}_{\bm{0}}$.
In particular, $\mathcal{X} \cap \mathbb{P}^{k-1} = \mathcal{Y}_{\bm{0}}$
and $\mathcal{X} = \mathbb{P}^{m-k} \ast \mathcal{Y}_{\bm{0}}$.
Since $\mathcal{X} \cap \mathbb{P}^{k-1} \supset \mathcal{Y}$
by the assumption
and $\deg \mathcal{X} = \deg \mathcal{Y}$,
we have
$\mathcal{Y}_{\bm{0}}= \mathcal{Y}$.
This completes the proof.
\qed

\section{A bound involving Koen Thas' invariant}
In \cite{tha2010}, Koen Thas defined an invariant of a hypersurface 
(see, Definition~\ref{KoenThas} below)
and obtained a bound for $N_q(\mathcal{X})$'s,
which involved the invariant.
We now give a simpler bound than his.
A comparison his bound and ours will give in the last section.

\begin{definition}\label{KoenThas}
Let $\mathcal{X}$ be a hypersurface defied over $\mathbb{F}_q$
in $\mathbb{P}^m$.
The maximum dimension of an $\mathbb{F}_q$-linear subspace
of $\mathbb{P}^m$ which is contained in $\mathcal{X}$ is denoted
by $k_{\mathcal{X}}$.
\end{definition}
By Lemma~\ref{linearinnonsingular}, if $\mathcal{X}$ is nonsingular
and $\deg \mathcal{X} \geq 2$, then
$k_{\mathcal{X}} \leq \lfloor \frac{m-1}{2}\rfloor.$

\begin{theorem}\label{theoremGEB}
Let $k$ be a nonnegative integer with $k \leq m-1$.
Let $\mathcal{X}$ be a hypersurface of degree $d$ over $\mathbb{F}_q$
in $\mathbb{P}^m$.
If $k_{\mathcal{X}} \leq k$, then
\begin{equation}\label{GEB}
N_q(\mathcal{X}) \leq \theta_q(m-k-1)\cdot q^k(d-1) + \theta_q(k).
\end{equation}
Furthermore, if $d\leq q$, the following conditions are equivalent:
\begin{enumerate}[{\rm (a)}]
 \item Equality holds in {\rm (\ref{GEB})};
 \item $k_{\mathcal{X}} = k$, and for any $\mathbb{F}_q$-linear subspace 
 $\mathcal{L}_1$ of dimension $k$ with $\mathcal{L}_1 \subset \mathcal{X}$
 and any $\mathcal{M}$ of dimension $k+1$
 with $\mathcal{L}_1 \subset \mathcal{M}$,
 \[(\ast)
   \left\{
  \begin{array}{l}
   \text{there are distinct } \mathbb{F}_q \text{-linear subspaces }
   \mathcal{L}_2, \dots , \mathcal{L}_d \text{ such that }\\
    \mathcal{M} \cap \mathcal{X} = \mathcal{L}_1 \cup \dots \cup \mathcal{L}_d
    \text{ and } \cap_{i=1}^d \mathcal{L}_i \text{ is of dimension }k-1.
  \end{array}
  \right.
  \]
 \item $k_{\mathcal{X}} = k$, and there is an $\mathbb{F}_q$-linear subspace 
 $\mathcal{L}_1$ of dimension $k$ with $\mathcal{L}_1 \subset \mathcal{X}$
 such that for any $\mathbb{F}_q$-linear subspace $\mathcal{M}$
 of dimension $k+1$ with $\mathcal{L}_1 \subset \mathcal{M}$,
 the condition $(\ast)$ is fulfilled.
\end{enumerate}
\end{theorem}
\proof
Put $\Phi(k,d)= \theta_q(m-k-1)\cdot q^k(d-1) + \theta_q(k)$.

{\em Step~1}. 
If $d \leq q+1$, then $\Phi(k+1, d) \geq \Phi(k, d).$
More precisely, if $d=q+1$, then $\Phi(k, q+1)=\theta_q(m)$ for any $k$;
and if $d \leq q$, then $\Phi(k+1, d) > \Phi(k, d).$

Actually, 
\begin{multline*}
\Phi(k+1, d) - \Phi(k, d) \\
        =\Bigl(\theta_q(m-(k+1)-1)\cdot q
                           - \theta_q(m-k-1) \Bigr) q^k(d-1)
                            + \theta_q(k+1) - \theta_q(k) \\
                         =-q^k(d-1) + q^{k+1} = q^k((q+1)-d), 
\end{multline*}
which is nonnegative if $d \leq q+1$, and
positive if $d < q+1$.
It is obvious that $\Phi(k, q+1)=\theta_q(m)$.

{\em Step~2}.
From Step~1, it is enough to show this theorem under
the assumption $k_{\mathcal{X}}=k$.
Choose any  $\mathbb{F}_q$-linear subspace
$\mathcal{L}_1$ of dimension $k$ with $\mathcal{L}_1 \subset \mathcal{X}$.
Let $\mathbb{G}$ be the set of $(k+1)$-dimensional $\mathbb{F}_q$-linear
subspaces containing $\mathcal{L}_1$.
Each point $P$ of $\mathcal{X} \setminus \mathcal{L}_1$
is contained in one and only one $(k+1)$-dimensional $\mathbb{F}_q$-linear
subspaces $\mathcal{M} \in \mathbb{G}$,
explicitly $\mathcal{M} = \langle \mathcal{L}_1 , P \rangle$.
Here $\langle \mathcal{L}_1 , P \rangle$ denotes the linear subspace
spanned by $\mathcal{L}_1$ and $P$.
Hence
\[
N_q(\mathcal{X}) = \sum_{\mathcal{M} \in \mathbb{G}}
     |(\mathcal{M}\cap\mathcal{X})(\mathbb{F}_q) \setminus 
                          \mathcal{L}_1(\mathbb{F}_q) |
       + N_q(\mathcal{L}_1).
\]
Applying the lemma of Segre-Serre-S{\o}rensen (\ref{theoremSSS})
for $\mathcal{M}\cap\mathcal{X} \subset \mathcal{M}= \mathbb{P}^{k+1}$,
\[
N_q(\mathcal{M}\cap\mathcal{X})
 \leq dq^k + \theta_q(k-1)
\]
and when $d \leq q$
equality holds if and only if the condition $(\ast)$ is satisfied.
On the other hand, $\mathbb{G}$ forms the set of  $\mathbb{F}_q$-points of
projective space $\mathbb{P}^{n-k-1}$.
Hence $|\mathbb{G}|= \theta_q(m-k-1)$ and
\begin{align*}
N_q(\mathcal{X}) \leq 
    & \theta_q(m-k-1)\cdot (dq^k + \theta_q(k-1) -\theta_q(k)) + \theta_q(k)\\
    =& \theta_q(m-k-1)\cdot q^k (d-1) + \theta_q(k)
\end{align*}
and when $d \leq q$
equality holds if and only if the condition $(\ast)$ is satisfied
for all $\mathcal{M} \in \mathbb{G}$.
This completes the proof.
\qed

\begin{remark}
If a hypersurface $\mathcal{X} \subset \mathbb{P}^{m}$
has no $\mathbb{F}_q$-hyperplane components,
then $k_{\mathcal{X}} \leq m-2$.
In this case, the bound (\ref{GEB}) is just the elementary bound
which we showed in \cite{hom-kim2013b}.
\end{remark}

\begin{corollary}\label{boundfornonsingular}
Let $\mathcal{X}$ be a nonsingular hypersurface of degree $d \geq 2$
of $\mathbb{P}^m$ over $\mathbb{F}_q$.
\begin{enumerate}[{\rm (i)}]
 \item If $m$ is odd, then
 \[
 N_q(\mathcal{X}) \leq \theta_q(\frac{m-1}{2})\cdot ((d-1)q^{\frac{m-1}{2}}+1).
 \]
 \item If $m$ is even, then
 \[
 N_q(\mathcal{X}) \leq \theta_q(\frac{m}{2})q^{\frac{m}{2}-1}(d-1)
                   + \theta_q(\frac{m}{2}-1).
 \]
\end{enumerate}
\end{corollary}
\proof
If $\mathcal{X}$ is nonsingular,
$k_{\mathcal{X}} \leq \lfloor \frac{m-1}{2} \rfloor$
by Lemma~\ref{linearinnonsingular}.
\qed

\section{Classification (the first step)}
By Lemma~\ref{linearinnonsingular},
in order to show the main theorem (Theorem~\ref{maintheorem}),
it is enough to prove the following theorem.
\begin{theorem}\label{actualtheorem}
Let $n$ be an odd integer at least $3$,
and $X$ a hypersurface of degree $d$ of $\mathbb{P}^n$
over $\mathbb{F}_q$.
If $k_X \leq \frac{n-1}{2}$,
then
\begin{equation}\label{bound}
N_q(X) \leq \theta_q(\frac{n-1}{2})\cdot ((d-1)q^{\frac{n-1}{2}}+1).
\end{equation}
Furthermore equality holds in {\rm (\ref{bound})}
if and only if
$X$ is one of the hypersurfaces in the list described
in Theorem~{\rm \ref{maintheorem}}.
\end{theorem}
The first part of this theorem has been already observed
in Corollary~\ref{boundfornonsingular}.

First we get rid of the cases $d=2$ and $d=q+1$.

\begin{proposition}\label{dis2}
Let $n$ be an odd integer at least $3$,
and $X$ a quadratic hypersurface of $\mathbb{P}^n$
over $\mathbb{F}_q$.
If $k_X \leq \frac{n-1}{2}$
and $N_q(X) = \theta_q(\frac{n-1}{2})(q^{\frac{n-1}{2}}+1),$
then $X$ is the nonsingular hyperbolic quadric,
that is, $X$ is projectively equivalent over $\mathbb{F}_q$
to the hypersurface
\[
\sum_{i=0}^{\frac{n-1}{2}} X_{2i} X_{2i+1} = 0.
\]
\end{proposition}
\proof
For a general theory of quadrics over a finite field,
consult \cite[Chapter~5]{hir1979}.
Since $k_X \leq \frac{n-1}{2}<n-1$,
the quadric does not split into two hyperplanes over
$\mathbb{F}_q$, that is, $X$ is irreducible over $\mathbb{F}_q$.
If $X$ is not absolutely irreducible,
then $X = H \cup H^{(q)}$ and
$X(\mathbb{F}_q) =  (H \cap H^{(q)})(\mathbb{F}_q)$,
where $H$ is a hyperplane over $\mathbb{F}_{q^2}$
and $H^{(q)}$ is the $q$-Frobenius conjugate of $H$.
This is a contradiction
because
$N_q(X) = \theta_q(\frac{n-1}{2})(q^{\frac{n-1}{2}}+1) =\theta_q(n-1) + q^{\frac{n-1}{2}}$
and $N_q(H \cap H^{(q)}) = \theta_q(n-2)$.
Therefore $X$ is absolutely irreducible,
and the possibilities of $X$ are as follows:
\begin{enumerate}[(i)]
 \item if $X$ is nonsingular, then $X$ is projectively equivalent
 over $\mathbb{F}_q$ to either
 \begin{align*}
  \mathcal{H}_n :& \sum_{i=0}^{\frac{n-1}{2}} X_{2i}X_{2i+1} =0 ;
     \text{ or}\\
  \mathcal{E}_n:& f(X_0, X_1) + \sum_{i=1}^{\frac{n-1}{2}} X_{2i}X_{2i+1} =0,
 \end{align*}
  where $f(X_0, X_1)$ is an irreducible quadratic polynomial 
  over  $\mathbb{F}_q$.
 \item if $X$ is a cone over a nonsingular quadric, then
 $X$ is projectively equivalent
 over $\mathbb{F}_q$ to either
 \begin{align*}
   \mathbb{P}^{n-2s-1} \ast \mathcal{P}_{2s}:& 
      X_0^2 + \sum_{i=1}^{s} X_{2i-1}X_{2i}=0
      \text{ with } s \leq \frac{n-1}{2}; \text{ or}\\
   \mathbb{P}^{n-2s} \ast \mathcal{H}_{2s-1}:& 
      \sum_{i=0}^{s-1} X_{2i}X_{2i+1}=0
      \text{ with } s \leq \frac{n-1}{2}; \text{ or}\\
  \mathbb{P}^{n-2s} \ast \mathcal{E}_{2s-1}:&
       f(X_0, X_1) + \sum_{i=1}^{s-1} X_{2i}X_{2i+1} =0
       \text{ with } s \leq \frac{n-1}{2}.
 \end{align*}
\end{enumerate}
If $X$ is one of the following quadrics:
  \begin{align*}
   \mathbb{P}^{n-2s-1} \ast \mathcal{P}_{2s}:& \text{ with }
              s \leq \frac{n-1}{2} -1 \text{ or} \\
    \mathbb{P}^{n-2s} \ast \mathcal{H}_{2s-1}:& \text{ with }
             s \leq \frac{n-1}{2}  \text{ or} \\
    \mathbb{P}^{n-2s} \ast \mathcal{E}_{2s-1}:& \text{ with }
            s \leq \frac{n-1}{2} -1,
  \end{align*}
then, $k_X > \frac{n-1}{2}$.
Actually, $\mathbb{P}^{n-2s-1} \ast \mathcal{P}_{2s}$
contains the $\mathbb{F}_q$-linear subspace
$X_0 = X_2 = X_4 =  \dots = X_{2s}=0$,
which is of dimension $n- (s+1)$, bigger than $\frac{n-1}{2}$
if $s \leq \frac{n-1}{2} -1$.
$\mathbb{P}^{n-2s} \ast \mathcal{H}_{2s-1}$ contains
$X_0 = X_2 = X_4 = \dots = X_{2(s-1)}=0$,
which is of dimension $n- s$, bigger than $\frac{n-1}{2}$
if $s \leq \frac{n-1}{2}$.
$\mathbb{P}^{n-2s} \ast \mathcal{E}_{2s-1}$ contains
$X_0 = X_1 = X_2 = X_4 = \dots = X_{2(s-1)}=0$,
which is of dimension $n- (s+1)$, bigger than $\frac{n-1}{2}$
if $s \leq \frac{n-1}{2} -1$.

So the remaining possibilities are either
$\mathcal{H}_n$ or $\mathcal{E}_n$ or $\mathbb{P}^0 \ast \mathcal{P}_{n-1}$
or $\mathbb{P}^1 \ast \mathcal{E}_{n-2}$.
Since
 \begin{align*}
    N_q(\mathcal{H}_n) &= \theta_q(\frac{n-1}{2})( q^{\frac{n-1}{2}} +1) 
                        = \theta_q(n-1) + q^{\frac{n-1}{2}}\\
    N_q(\mathcal{E}_n) &= \theta_q(\frac{n-3}{2})( q^{\frac{n+1}{2}} +1) 
                         = \theta_q(n-1) - q^{\frac{n-1}{2}}\\
    N_q(\mathbb{P}^0 \ast \mathcal{P}_{n-1})&=
          N_q(\mathcal{P}_{n-1})q+1 = \theta_q(n-2)q +1 
                             = \theta_q(n-1)\\
    N_q(\mathbb{P}^1 \ast \mathcal{E}_{n-2}) &=
          N_q(\mathcal{E}_{n-2})q^2 + \theta_q(1)
                  = \theta_q(n-1) - q^{\frac{n+1}{2}},
 \end{align*}
$X$ must be projectively equivalent over $\mathbb{F}_q$
to $\mathcal{H}_n$.
\qed

\begin{proposition}
Let $n$ be an odd integer at least $3$,
and $X$ a hypersurface of degree $q+1$ of $\mathbb{P}^n$
over $\mathbb{F}_q$.
If $k_X \leq \frac{n-1}{2}$ and
$N_q(X) = \theta_q(\frac{n-1}{2})\dot (q^{\frac{n-1}{2}}+1),$
then $X$ is projectively equivalent over $\mathbb{F}_q$ to
the hypersurface
\[
\sum_{i=0}^{\frac{n-1}{2}} (X_{2i}^q X_{2i+1} - X_{2i}X_{2i+1}^q )=0.
\]
\end{proposition}
\proof
Since $\theta_q(\frac{n-1}{2})\dot (q^{\frac{n-1}{2}}+1)= \theta_q(n),$
$X(\mathbb{F}_q) = \mathbb{P}^n(\mathbb{F}_q).$
Hence the ideal of $X$ is generated by
$\{ X_i^q X_j - X_i X_j^q \mid i < j \}$.
Therefore, there is a $q$-alternating matrix $A$ over $\mathbb{F}_q$
such that $X$ is given by the equation
\[
(X_0^q, \dots , X_n^q) A 
                  \left(
                     \begin{array}{c}
                       X_0 \\
                       \vdots \\
                       X_n
                     \end{array}
                  \right)
                                  =0.
\]
By the standard theory of alternating matrix over $\mathbb{F}_q$,
we can choose new coordinates $X_0, \dots , X_n$
of $\mathbb{P}^n$ over $\mathbb{F}_q$ so that
$A$ is of the form
\[
\left(
   \begin{array}{rccrcc}
   0 & 1&       &  &    &   \\
   -1& 0&       &  &    &   \\
     &  & \ddots&  &    &   \\
     &  &       &0 & 1 &    \\
     &  &       &-1& 0 &    \\
     &  &       &  &   & \Huge{\rm O}
   \end{array}
\right),
\]
that is,
$X$ is defined by 
$
\sum_{i=0}^s (X_{2i}^q X_{2i+1} - X_{2i}X_{2i+1}^q )=0
$
with $ s \leq \frac{n-1}{2}.$
Obviously,
$
\{ X_0 = X_2 = \dots = X_{2s} = 0 \} \subset X,
$
and this $\mathbb{F}_q$-linear subspace is
of dimension $n-(s+1)$.
Since $k_X \leq \frac{n-1}{2}$,
we have $s=\frac{n-1}{2}$.
\qed

\section{Classification (continuation)}
To complete the proof of Theorem~\ref{actualtheorem},
we clarify the necessary set-up.
In the previous section, two cases $d=2 \text{ and } q+1$
were  already handled. 
\begin{setup}\label{setup1}
Let $n$ be an odd integer at least $3$,
and $X$ a hypersurface of $\mathbb{P}^n$
over $\mathbb{F}_q$.
Suppose that the degree $d$ of $X$ is in the range $2<d \leq q$,
$k_X = \frac{n-1}{2}$ and 
\begin{equation}\label{eqsetup1}
N_q(X) = \theta_q( \frac{n-1}{2})\cdot ( (d-1)q^{ \frac{n-1}{2}} +1).
\end{equation}
\end{setup}
Note that initially the condition $k_X  \leq \frac{n-1}{2}$
was supposed in Theorem~\ref{actualtheorem}, however, since we may assume
that $d \leq q$ at this stage, the condition $k_X = \frac{n-1}{2}$
holds by Theorem~\ref{theoremGEB}.

\begin{notation}
The set of $\mathbb{F}_q$-linear subspaces of dimension $u$ in $\mathbb{P}^n$
is denoted by $G(u, \mathbb{P}^n)(\mathbb{F}_q)$.
\end{notation}

\begin{definition}
For $X$ in Set-up~\ref{setup1},
$M \in  G(\frac{n+1}{2}, \mathbb{P}^n)(\mathbb{F}_q)$
is said to be of type S (for $X$) if
\[
M \cap X = L_1 \cup \dots \cup L_d,
\]
where $L_1 , \dots , L_d \in  G(\frac{n-1}{2}, \mathbb{P}^n)(\mathbb{F}_q)$
and $\cap_{i=1}^d L_i \in G(\frac{n-3}{2}, \mathbb{P}^n)(\mathbb{F}_q)$.
This $\frac{n-3}{2}$-dimensional linear subspace is denoted by $\Lambda_M$.
\end{definition}

The number of $\mathbb{F}_q$-points of
$M\cap X$ above is given by:
\begin{lemma}\label{numberofMX}
\[
  |(M \cap X)(\mathbb{F}_q)| =
      d q^{\frac{n-1}{2}} + \theta_q(\frac{n-3}{2}),
\]
\end{lemma}
\proof
This is a direct consequence of Lemma~\ref{dL}.
\qed

\begin{lemma}\label{existenceL}
Let $M \in  G(\frac{n+1}{2}, \mathbb{P}^n)(\mathbb{F}_q)$.
Then there is a linear space
$L_1 \in  G(\frac{n-1}{2}, \mathbb{P}^n)(\mathbb{F}_q)$
with $L_1 \subset X$ such that $L_1 \subset M$ if and only if 
$M$ is of type S.
\end{lemma}
\proof
The {\em if} part is obvious by definition.
The {\em only if} part comes from
Theorem~\ref{theoremGEB}, (a) $\Rightarrow$ (b).
\qed

\begin{remark}\label{singularlocusofMcapX}
When $M \in  G(\frac{n+1}{2}, \mathbb{P}^n)(\mathbb{F}_q)$
is of type S, then $\Sing (M\cap X) = \Lambda_M$
by Lemma~\ref{singularlocus}.
\end{remark}

We need further notation:
\begin{notation}
   \begin{itemize}
      \item $\mathbb{L}:=\{
        L \in G( \frac{n-1}{2}, \mathbb{P}^n )(\mathbb{F}_q) \mid L \subset X
                         \}.$
      \item For $P \in X(\mathbb{F}_q),$
       $\mathbb{L}(P) :=\{ L \in \mathbb{L} \mid L \ni P \}.$
   \end{itemize}
\end{notation}

\begin{lemma}\label{LPisnonempty}
 For any $P \in X(\mathbb{F}_q)$, $\mathbb{L}(P)\neq \emptyset .$
\end{lemma}
\proof
By Theorem~\ref{theoremGEB},
$\mathbb{L} \neq \emptyset .$
Choose $L_1 \in \mathbb{L}$.
Then either $P \in L_1$ or $P \not\in L_1$.
If the latter case occurs, then
$M =\langle L_1, P\rangle \in G( \frac{n+1}{2}, \mathbb{P}^n )(\mathbb{F}_q)$.
Then
$M\cap X = L_1 \cup \dots \cup L_d$
by (\ref{existenceL}).
Hence $P \in L_i$ for some $i$, that is, $L_i \in \mathbb{L}(P).$
\qed

\begin{lemma}\label{nonsingularpoint}
 Let $L \in \mathbb{L}.$
 If $P \in X(\mathbb{F}_q) \setminus L,$
 then $P$ is a nonsingular point of $X$.
\end{lemma}
\proof
Let $M=\langle L, P\rangle 
      \in G( \frac{n+1}{2}, \mathbb{P}^n )(\mathbb{F}_q),$
which is of type S by (\ref{existenceL}).
Since $\Sing (M\cap X) = \Lambda_M \subset L$ by (\ref{singularlocusofMcapX}),
$P$ is a nonsingular point of $M \cap X$.
Hence so is $P$ in $X$ by (\ref{nonsingularity_ancestor}).
\qed

\begin{proposition}\label{presymmetric}
  Let $P_0$ be an $\mathbb{F}_q$-point of $X$.
  Suppose $P_0$ is a nonsingular point of $X$.
  \begin{enumerate}[{\rm (i)}]
   \item If $L_1 \in \mathbb{L}(P_0),$
      then $L_1 \subset T_{P_0}X$,
      where $T_{P_0}X$ is the embedded tangent hyperplane to $X$ at $P_0$.
   \item Let $L_1 \in \mathbb{L}(P_0)$, and $M$ of type S containg $L_1$.
   If $M \subset T_{P_0}X$, then $P_0 \in \Lambda_M$.
   \item If $M$ is of type S and $\Lambda_M \ni P_0,$
      then $M \subset T_{P_0}X.$
  \end{enumerate}
\end{proposition}
\proof
(i) Since $P_0 \in L_1 \subset X$,
we have
$T_{P_0} L_1 =L_1$ (because $L_1$ itself is linear) and
$T_{P_0}L_1 \subset T_{P_0}X$.
Hence
$L_1 \subset T_{P_0}X$.

(ii) Since $M$ is of type S containing $L_1$,
there are $L_2, \dots ,  L_d  \in \mathbb{L}$
such that
$
M \cap X = L_1 \cup L_2 \cup \dots \cup L_d .
$
Since $P_0$ is a singular point of $T_{P_0}X \cap X$
which is  a hypersurface
of $\mathbb{P}^{n-1}=T_{P_0}X$,
it is also singular point of $(T_{P_0}X \cap X) \cap M$
by (\ref{nonsingularity_ancestor}).
Since $(T_{P_0}X \cap X) \cap M = X \cap M$
(because the assumption $M \subset T_{P_0}X$),
$P_0 \in \Sing (X \cap M) = \Lambda_M$.

(iii)
There are $\mathbb{F}_q$-linear subspaces
 $L_1, \dots ,  L_d  \in \mathbb{L}$
 such that 
$
M \cap X = L_1 \cup L_2 \cup \dots \cup L_d 
$
with $\Lambda_M = \cap_{i=1}^d L_i$.
Hence $P_0 \in L_i$ for any $i = 1, \dots , d.$
Since $L_i \subset T_{P_0}$ by (i) and
$\langle L_1, \dots ,  L_d  \rangle = M,$
we have $M \subset T_{P_0}X$
\qed

\begin{corollary}\label{symmetric}
 Let $P_0, P_1 \in X(\mathbb{F}_q)$ be two distinct nonsingular points
 of $X$. Then
 $T_{P_0}X \ni P_1$ if and only if $T_{P_1}X \ni P_0$.
\end{corollary}
\proof
Suppose the condition $T_{P_0}X \ni P_1$.
We can find an $\mathbb{F}_q$-space $L_1 \in \mathbb{L}(P_0)$ by
Lemma~\ref{LPisnonempty}.
When $P_1 \in L_1$, $L_1 = T_{P_1}L_1 \subset T_{P_1}X$.
Since $P_0 \in L_1$, we have $P_0 \in T_{P_1}X$.
When $P_1 \not\in L_1$, let $M: =\langle L_1 , P_1 \rangle$.
Since $L_1= T_{P_0}L_1 \subset T_{P_0}X$ and
$P_1 \in T_{P_0}X$ by the assumption, we have $M \subset T_{P_0}X$.
Hence $P_0 \in \Lambda_M$ by (ii) of Proposition~\ref{presymmetric}.
Since
$M \cap X = L_1 \cup L_2^{(M)} \cup \dots \cup L_d^{(M)}$
where $L_i^{(M)} \in \mathbb{L}$ ($i=2, \dots , d$),
there is an $L_i^{(M)}$ which contains $P_1$.
Hence $L_i^{(M)} \subset T_{P_1}X$.
On the other hand,
since $P_0 \in \Lambda_M \subset L_i^{(M)}$,
we can conclude that $P_0 \in T_{P_1}X$.
\qed

\begin{setup}\label{setup2}
We keep Set-up~\ref{setup1}.
Additionally, fix a nonsingular point $P_0 \in X(\mathbb{F}_q)$
(the existence of such a point has been guaranteed
by Lemma~\ref{nonsingularpoint} and (\ref{eqsetup1})),
and also $L_1 \in \mathbb{L}(P_0)$.
Let $Y$ be the hypersurface
$X \cap T_{P_0}X$ in $T_{P_0}X = \mathbb{P}^{n-1}$,
which is also defined over $\mathbb{F}_q$ and of degree $d$.
\end{setup}

\begin{lemma}\label{numberofY}
\[
N_q(Y) = \theta_q(\frac{n-3}{2})\cdot 
            q^{\frac{n-1}{2}} (d-1) + \theta_q(\frac{n-1}{2}).
\]
\end{lemma}
\proof
Let
\[
\mathbb{G} = \{
M \in G(\frac{n+1}{2}, \mathbb{P}^n)(\mathbb{F}_q)
       \mid L_1 \subset M \subset T_{P_0}X
\}.
\]
Then $\mathbb{G}$ forms a finite projective space
$\mathbb{P}^{\frac{n-3}{2}}(\mathbb{F}_q)$.
Obviously, 
$Y(\mathbb{F}_q) = \cup_{M \in \mathbb{G}}(M \cap X)(\mathbb{F}_q)$
and $M \cap M' =L_1$ if $M$ and $M'$ are distinct elements of $\mathbb{G}$.
Hence
\begin{align*}
|Y(\mathbb{F}_q)| &= \sum_{M \in \mathbb{G}}
         \left(  |(M \cap X)(\mathbb{F}_q)|-|L_1(\mathbb{F}_q)| \right)
                                                     +|L_1(\mathbb{F}_q)| \\
    &=\theta_q(\frac{n-3}{2})q^{\frac{n-1}{2}} (d-1) + \theta_q(\frac{n-1}{2}),
\end{align*}
where the last equality comes from Lemma~\ref{numberofMX}.
\qed

\begin{setup}\label{setup3}
We keep Set-ups \ref{setup1} and \ref{setup2}.
Furthermore, suppose $n \geq 5$.
Take an $\mathbb{F}_q$-hyperplane $H \subset \mathbb{P}^n$
so that $H \not\ni P_0$.
Then $T_{P_0}X \cap H$ is a linear subspace defined over $\mathbb{F}_q$
of codimension $2$ in $\mathbb{P}^n$.
Let $Z$ be the hypersurface 
\[
Y \cap \left( T_{P_0}X \cap H \right)
\text{ in }
T_{P_0}X \cap H = \mathbb{P}^{n-2},
\]
which is also defined over $\mathbb{F}_q$
and of degree $d$.
Note that since $Y \subset T_{P_0}X$,
$Z$ is just a cutout of $Y$ by $H$, that is,
$Z=Y \cap H.$
\end{setup}

\begin{lemma}\label{numberofZ}
\[
N_q(Z) = \theta_q(\frac{n-3}{2}) \cdot \Bigl((d-1)q^{\frac{n-3}{2}} +1 \Bigr).
\]
\end{lemma}
\proof
Since
\[
Y(\mathbb{F}_q) = 
  \bigcup_{M \in \mathbb{G}} 
     (L_1 \cup L_2^{(M)} \cup \dots \cup L_d^{(M)})(\mathbb{F}_q),
\]
we have
\[
Z(\mathbb{F}_q) = \bigcup_{M \in \mathbb{G}} 
 \Bigl((L_1\cap H) \cup (L_2^{(M)}\cap H) \cup 
                        \dots \cup (L_d^{(M)}\cap H) \Bigr)(\mathbb{F}_q).
\]
Since $(M \cap H) \cap (M' \cap H)= L_1 \cap H$
if $M$ and $M'$ are distinct elements of $\mathbb{G}$,
\begin{equation}\label{countingZ}
  \begin{split}
  | Z(\mathbb{F}_q)| = &\sum_{M \in \mathbb{G}} 
  \Bigl(
 |((L_1\cap H) \cup (L_2^{(M)}\cap H) \cup 
                   \dots \cup (L_d^{(M)}\cap H))(\mathbb{F}_q)|    \\
                         &  -|(L_1\cap H)(\mathbb{F}_q)| \Bigr) 
                        + |(L_1\cap H)(\mathbb{F}_q)|.
  \end{split}
\end{equation}
For each $M \in \mathbb{G}$,
since $\Lambda_M \ni P_0$ (\ref{presymmetric}, ii) but $H \not\ni P_0$,
\[
\dim \, L_1\cap H = \dim \, L_2^{(M)}\cap H =
                      \dots = \dim\,  L_d^{(M)}\cap H = \frac{n-3}{2},
\]
and
\[
(L_1\cap H) \cap (L_2^{(M)}\cap H) \cap 
                           \dots \cap (L_d^{(M)}\cap H)
                     = \Lambda_M \cap H = \mathbb{P}^{\frac{n-5}{2}}.
\]
Hence
\begin{equation}\label{countingLH}
  |((L_1\cap H) \cup (L_2^{(M)}\cap H) \cup 
                   \dots \cup (L_d^{(M)}\cap H))(\mathbb{F}_q)|  
                               =dq^{\frac{n-3}{2}} + \theta_q(\frac{n-5}{2}) 
\end{equation}
by Lemma~\ref{dL}.
Therefore, by (\ref{countingZ}) and (\ref{countingLH})
\begin{align*}
N_q(Z) =& \theta_q(\frac{n-3}{2}) 
         \Bigl( d q^{\frac{n-3}{2}} +  \theta_q(\frac{n-5}{2}) 
           -  \theta_q(\frac{n-3}{2})\Bigr) 
                 +\theta_q(\frac{n-3}{2}) \\
          = & \theta_q(\frac{n-3}{2}) \Bigl( (d-1)q^{\frac{n-3}{2}} +1 \Bigr).
           \text{ \ \ \     \qed}
\end{align*}

\begin{lemma}\label{kZ}
\[
k_Z = \frac{n-3}{2}.
\]
\end{lemma}
\proof
Since $L_1 \cap H \subset Z$,
$k_Z \geq \frac{n-3}{2}.$

Suppose there is an $\frac{n-1}{2}$-dimensional $\mathbb{F}_q$-linear
space $L_0$ which is contained in $Z\subset X$.
Then for each $Q \in Z(\mathbb{F}_q) \setminus L_0$,
$M:= \langle L_0, Q \rangle$ is of type S for $X$,
and is contained in $T_{P_0}X \cap H = \mathbb{P}^{n-2}$
(because $L_0 \subset Z$
and $Q \in Z$).
Let
\begin{align*}
 \mathbb{G}'& := \{
                  M \in G(\frac{n+1}{2}, \mathbb{P}^{n})(\mathbb{F}_q)
                    \mid L_0 \subset M \subset \mathbb{P}^{n-2}
                        = T_{P_0}X \cap H
                \} \\
            &= \{
                 M \in G(\frac{n+1}{2}, \mathbb{P}^{n-2})(\mathbb{F}_q)
                   \mid L_0 \subset M 
               \}\\
            &= \mathbb{P}^{\frac{n-5}{2}}(\mathbb{F}_q).
\end{align*}
Since
 \begin{enumerate}[(i)]
    \item $Z(\mathbb{F}_q) = \cup_{M \in \mathbb{G}'}
               (M \cap X)(\mathbb{F}_q),$
    \item $M \cap M' = L_0$ for distinct elements 
                   $M, M' \in \mathbb{G}'$  and
    \item $|(M \cap X)(\mathbb{F}_q)| =
      d q^{\frac{n-1}{2}} + \theta_q(\frac{n-3}{2})$
 by Lemma~\ref{numberofMX},
 \end{enumerate}
we can compute the number of $Z(\mathbb{F}_q)$ as
\begin{align}
 Z(\mathbb{F}_q) = & \theta_q(\frac{n-5}{2}) 
           \Bigl(
            d q^{\frac{n-1}{2}} +\theta_q(\frac{n-3}{2})
                   - |L_0(\mathbb{F}_q)|
           \Bigr) + |L_0(\mathbb{F}_q)|  \notag \\
                 & = \theta_q(\frac{n-5}{2})(d-1)q^{\frac{n-1}{2}}
                    + \theta_q(\frac{n-1}{2}).
                    \label{thisnumber}
\end{align}
Compare this number (\ref{thisnumber})
with that computed in Lemma~\ref{numberofZ}.
Namely,
\begin{align*}
 \Bigl(
 \theta_q(\frac{n-5}{2})(d-1)q^{\frac{n-1}{2}}
                    + \theta_q(\frac{n-1}{2})
 \Bigr) -&
 \Bigl(
  \theta_q(\frac{n-3}{2})
    \bigl(
     (d-1)q^{\frac{n-3}{2}} +1
    \bigr)
 \Bigr) \\
 &= q^{\frac{n-3}{2}}(q+1 -d),
\end{align*}
which is a contradiction because $d \leq q$.
Therefore $k_Z = \frac{n-3}{2}$.
\qed

\begin{theorem}
Under Set-up~{\rm \ref{setup1}},
$q$ is square and $d=\sqrt{q} +1$.
\end{theorem}
\proof
When $n=3$, we already know that
the conclusion is true (Theorem~\ref{nequal3}).
By Lemmas~\ref{numberofZ} and \ref{kZ},
the induction on odd $n$ works well.
\qed
\section{Classification for $d= \sqrt{q} +1$}
The remaining part of the classification is
to determine the structure of $X$ under Set-up~\ref{setup1}
when $d= \sqrt{q}+1$.
Of course, throughout this section, $q$ is supposed to be square.

When $n=3$, we already know the surface $X$ is a nonsingular Hermitian
surface \cite{hom-kim2015online}.
So we suppose that $n \geq 5$ as we did after Set-up~\ref{setup3}.
We keep the situation described in Set-ups~\ref{setup1} and \ref{setup2}.

\begin{lemma}\label{canchoose}
The set $X(\mathbb{F}_q) \setminus T_{P_0}X$ is nonempty, and
each point of this set is a nonsingular points of $X$.
\end{lemma}
\proof
Note that $X(\mathbb{F}_q) \setminus T_{P_0}X = 
                  X(\mathbb{F}_q) \setminus Y(\mathbb{F}_q)$
       because $Y = X \cap T_{P_0}X$
        (see Set-up~\ref{setup2}).
By Set-up~\ref{setup1} and Lemma~\ref{numberofY},
\begin{align*}
 N_q(X) - N_q(Y) &=  \\
    & \theta_q(\frac{n-1}{2}) \bigl(
                                (d-1)q^{\frac{n-1}{2}}+1
                              \bigr)
      -\Bigl(\theta_q(\frac{n-3}{2})q^{\frac{n-1}{2}}(d-1)
                + \theta_q(\frac{n-1}{2}) \Bigr) \\
      & = (d-1)q^{n-1}  = q^{n-\frac{1}{2}}>0.
\end{align*}
Hence $X(\mathbb{F}_q) \setminus T_{P_0}X \neq \emptyset .$
Since $L_1 \in \mathbb{L}$ lies on $T_{P_0} X$
by Proposition ~\ref{presymmetric} (i),
any point of $X(\mathbb{F}_q) \setminus T_{P_0}X$
is nonsingular by Lemma~\ref{nonsingularpoint}.
\qed

\begin{proposition}\label{inductionstepZ}
Suppose $n$ is an odd integer with $n \geq 5$.
Let $X$ be a hypersurface of degree $\sqrt{q}+1$
in $\mathbb{P}^n$ over $\mathbb{F}_q$ with the conditions described
in Set-up~{\rm \ref{setup1}}.
Let $Q_0$ and $Q_1$ be points of $X(\mathbb{F}_q)$ that are nonsingular points
of $X$ with $T_{Q_1} \not\ni Q_0.$ {\rm(}Hence $T_{Q_0} \not\ni Q_1$ neither
by Corollary~{\rm \ref{symmetric}}.{\rm )}
Let $Y = X \cap T_{Q_0}X$ , $Y' = X \cap T_{Q_1}X$,
and
\[
Z = Y \cap T_{Q_1}X = Y' \cap T_{Q_0}X = X \cap T_{Q_0}X \cap T_{Q_1}X.
\]
Then 
$
N_q(Z) = \theta_q(\frac{n-3}{2}) \bigl(q^{\frac{n-2}{2}} +1 \bigr)
$
and
$
k_Z = \frac{n-3}{2}.
$
Furthermore, 
$ Y = Q_0 \ast Z$ in $T_{Q_0}X =\mathbb{P}^{n-1}$
and $ Y' = Q_1 \ast Z$ in $T_{Q_1}X =\mathbb{P}^{n-1}.$
\end{proposition}
\proof
Regard $Q_0$ as the point $P_0$ in Set-ups~\ref{setup2} and \ref{setup3},
and $T_{Q_1}X$ as the hyperplane $H$.
Then
$
N_q(Z) = \theta_q(\frac{n-3}{2}) \bigl(q^{\frac{n-2}{2}} +1 \bigr)
$
by Lemma~\ref{numberofZ} with the assumption $d= \sqrt{q}+1$,
and also 
$
k_Z = \frac{n-3}{2}.
$
by Lemma~\ref{kZ}.

Choose coordinates $X_1, \dots , X_n$
of $T_{Q_0}X = \mathbb{P}^{n-1}$
so that $Q_0 =(1, 0, \dots , 0)$ in $\mathbb{P}^{n-1}$
and $T_{Q_0}X \cap T_{Q_1}X = \{ X_1 = 0 \}$ in $\mathbb{P}^{n-1}$.
We want to apply the cone lemma (Proposition~\ref{conelemma})
to our situation, that is,
regard the hypersurface $Y$ of $\mathbb{P}^{n-1} = T_{Q_0}X$
as the hypersurface $\mathcal{X}$ of $\mathbb{P}^{m}$
in (\ref{conelemma}), 
$Z \subset \mathbb{P}^{n-2}= T_{Q_0}X \cap \{ X_1=0 \}$
as $\mathcal{Y} \subset \mathbb{P}^{k-1}$, 
and $Q_0 = \mathbb{P}^0$ as $\mathcal{L}=\mathbb{P}^{m-k}.$
So $m$ and $k$ in the cone lemma are both $n-1$ in the current situation.
The first condition in (\ref{conelemma})
can be paraphrased in our situation as
\[
N_q(Z) =
\theta_q(\frac{n-3}{2}) \bigl(
                          q^{\frac{n-2}{2}} + 1
                         \bigr)
             > \sqrt{q}q^{n-3} + \theta_q(n-4),
\]
and it is not hard to check this inequality holds.
The second condition in (\ref{conelemma})
obviously holds.
To check the last condition, let $R \in Z(\mathbb{F}_q).$
Choose $L_1 \in \mathbb{L}(Q_0)$,
and let $M= \langle  L_1 , R \rangle \subset T_{Q_0}X$
if $R \not\in L_1$.
Then 
$M \cap X = L_1^{(M)} \cup L_2^{(M)} \cup \dots \cup L_d^{(M)}
    \subset T_{Q_0}X,$
and $Q_0 \in \Lambda_M = \cap_{i=1}^{d}L_i^{(M)},$
where $L_1^{(M)} =L_1.$
Since there is an $L_i^{(M)}$
such that $R \in L_i^{(M)}$,
the line $\langle Q_0 , R \rangle$
is contained in $L_i^{(M)}$.
Since $L_i^{(M)} \subset T_{Q_0}X \cap X =Y$,
we can conclude that $(Q_0 \ast Z)(\mathbb{F}_q) \subset Y.$

Therefore, by the cone lemma, 
$ Y = Q_0 \ast Z$.
By the symmetry of the role of $Q_0$ and that of $Q_1$,
$ Y' = Q_1 \ast Z$ also holds.
\qed

\gyokan

We finally prove the following theorem which completes the proof
of Theorem~\ref{actualtheorem}.

\begin{theorem}
Suppose $n$ is an odd integer with $n \geq 3$.
Let $X$ be a hypersurface of degree $\sqrt{q} +1$
in $\mathbb{P}^n$ defined over $\mathbb{F}_q$.
If $k_X=\frac{n-1}{2}$ and
$
N_q(X) = \theta_q(\frac{n-1}{2}) \bigl(q^{\frac{n}{2}} +1 \bigr),
$
then $X$ is a nonsingular Hermitian hypersurface.
\end{theorem}
\proof
When $n=3$, this was already proved in \cite{hom-kim2015online}.
So we assume that $n \geq 5$.

First we choose a point $P_0 \in X(\mathbb{F}_q)$
which fits with Set-ups~\ref{setup1} and \ref{setup2}.
By Lemma~\ref{canchoose}, we can choose a point
$P_1 \in X(\mathbb{F}_q) \setminus T_{P_0}X$,
and it is a nonsingular point of $X$.
Hence $P_0 \not\in T_{P_1}X$ by Corollary~\ref{symmetric}.
Choose coordinates $X_0, X_1, \dots , X_n$ of $\mathbb{P}^n$
over $\mathbb{F}_q$
so that $P_0 = (0,1,0, \dots , 0)$,
$P_1 = (1,0, \dots , 0)$,
$T_{P_0}X= \{  X_0= 0\}$ and
$T_{P_1}X= \{  X_1= 0\}$.
Note that if one applies a linear transformation of type
\[
\begin{pmatrix}
1_2 & 0 \\
0   & A
\end{pmatrix}
\, \, 
(\,  A \in GL(n-1, \mathbb{F}_q) \,),
\]
to these coordinates,
it does not affect the coordinate representations of
$P_0$ and $P_1$,
and the equations of $T_{P_0}X$ and $T_{P_1}X$.

Let $Y= X \cap T_{P_0}X$,  $Y'= X \cap T_{P_1}X$
and $Z= X \cap T_{P_0}X \cap T_{P_1}X$.
Since $\mathbb{P}^{n-2} = T_{P_0}X \cap T_{P_1}X$
is defined by $X_0 = X_1=0$,
we can regard $X_2, \dots , X_n$ as coordinates of
$T_{P_0}X \cap T_{P_1}X$.
By Proposition~\ref{inductionstepZ},
we can apply the induction hypothesis to $Z$,
that is, $Z$ is a nonsingular Hermitian hypersurface
in $\mathbb{P}^{n-2} = T_{P_0}X \cap T_{P_1}X$.
Therefore, we may assume that
$Z$ is defined by
\begin{equation}\label{eqZ}
\sum_{i=1}^{\frac{n-1}{2}}
\Bigl(
X_{2i}^{\sqrt{q}}X_{2i+1} + X_{2i}X_{2i+1}^{\sqrt{q}}
\Bigr)
  = 0.
\end{equation}
Since $Y= P_0 \ast Z$ and
$Y'= P_1 \ast Z$
in $T_{P_0}X = \mathbb{P}^{n-1}$
and $T_{P_1}X = \mathbb{P}^{n-1}$ respectively,
the equation (\ref{eqZ}) is also
that for $Y$ with coordinates $X_0, X_2, \cdots , X_n$
and that for $Y'$ with coordinates $X_1, \cdots , X_n$
respectively.
Therefore $X$ is defined by
$F=0$ with
\begin{equation}\label{eqF}
F=X_0 X_1 G(X_0, \dots , X_n)
+ \sum_{i=1}^{\frac{n-1}{2}}
\Bigl(
X_{2i}^{\sqrt{q}}X_{2i+1} + X_{2i}X_{2i+1}^{\sqrt{q}}
\Bigr),
\end{equation}
where $G(X_0, \dots , X_n)$ is a homogeneous polynomial
of degree $\sqrt{q}-1$.
The partial derivations of $F$ are as follows:
\begin{equation}\label{pdF}
\begin{split}
  &\cfrac[l]{\partial F}{\partial X_0} 
         = X_1 G +
         X_0 X_1 \cfrac[l]{\partial G}{\partial X_0} \\
  &\cfrac[l]{\partial F}{\partial X_1} 
         = X_0 G +
         X_0 X_1 \cfrac[l]{\partial G}{\partial X_1} \\
  &\cfrac[l]{\partial F}{\partial X_{2i}} 
         = X_0 X_1 \cfrac[l]{\partial G}{\partial X_{2i}}
              + X_{2i+1}^{\sqrt{q}}
               \qquad (1 \leq i \leq \frac{n-1}{2}) \\
  &\cfrac[l]{\partial F}{\partial X_{2i+1}} 
         = X_0 X_1 \cfrac[l]{\partial G}{\partial X_{2i+1}}
              + X_{2i}^{\sqrt{q}}
               \qquad (1 \leq i \leq \frac{n-1}{2}).
\end{split}
\end{equation}
For each $i=1, 2, \dots , \frac{n-1}{2}$,
let
\begin{align*}
    P_{2i} &= (0, \dots , 0 , \overset{2i}{0},
                          \overset{2i+1}{1}, 0, \dots , 0)\\
    P_{2i+1} &= (0, \dots , 0 , \overset{2i}{1} ,
                          \overset{2i+1}{0}, 0, \dots , 0).
\end{align*}
Then these points are nonsingular points of $X$,
$T_{P_{2i}}X = \{X_{2i}=0\}$, and
$T_{P_{2i+1}}X = \{X_{2i+1}=0\}$ by (\ref{pdF}).
Apply Proposition~\ref{inductionstepZ} to $P_{2i}$
and $P_{2i+1}$.
Then $X \cap T_{P_{2i}}X \cap T_{P_{2i+1}}X$
is also a nonsingular Hermitian hypersurface
in $T_{P_{2i}}X \cap T_{P_{2i+1}}X = \mathbb{P}^{n-2}$
by the induction hypothesis.

Here we need a little more terminology:
for letters $X_0, \dots , X_n$
over $\mathbb{F}_q$, polynomials of type
\[
X_k^{\sqrt{q}+1} \qquad \text{or} \qquad 
   \lambda X_k^{\sqrt{q}} X_l + \lambda^{\sqrt{q}} X_k X_l^{\sqrt{q}}
   \quad (\lambda \in \mathbb{F}_q^{\times})
\]
are referred as Hermitian molecules.
An equation of a Hermitian hypersurface, by definition,
consists of an $\mathbb{F}_q$-linear combination of Hermitian molecules
(but the converse is not true).

Since
\begin{equation}\label{eqsubF}
F(X_0, \dots , X_{2i-1}, \overset{2i}{0}, \overset{2i+1}{0}, 
                   X_{2i+1}, \dots , X_n) = 0
\end{equation}
is an equation of the Hermitian hypersurface
$X \cap T_{P_{2i}}X \cap T_{P_{2i+1}}X$
in $\mathbb{P}^{n-2}$,
\[
X_0 X_1 G(X_0, \dots , X_{2i-1}, \overset{2i}{0}, \overset{2i+1}{0}, 
                   X_{2i+1}, \dots , X_n)
\]
consists of Hermitian molecules.
Hence
\begin{equation}\label{eqsubG}
G(X_0, \dots , X_{2i-1}, \overset{2i}{0}, \overset{2i+1}{0}, 
                   X_{2i+1}, \dots , X_n)
  = c \bigl( \lambda X_0^{\sqrt{q}} + \lambda^{\sqrt{q}} X_1^{\sqrt{q}} \bigr)
\end{equation}
for appropriate $\lambda \in \mathbb{F}_q$ and $c \in \overline{\mathbb{F}}_q$.
Since the equation (\ref{eqsubF})
defines a Hermitian hypersurface and
the polynomial
contains a pair of terms
$X_{2j}^{\sqrt{q}}X_{2j+1} + X_{2j}X_{2j+1}^{\sqrt{q}}$
for some $j \geq 1$, we know
$ c \in \mathbb{F}_{\sqrt{q}},$
that is, we may assume $c$ to be $1$ in (\ref{eqsubG}),
and also
\begin{equation}\label{firsttermsofF}
X_0X_1G(X_0, \dots , X_n)
 = X_0X_1 
  \bigl( \lambda X_0^{\sqrt{q}} + \lambda^{\sqrt{q}} X_1^{\sqrt{q}} \bigr)
 +H(X_0, \dots, X_n)
\end{equation}
with 
\begin{equation}\label{conditionH}
H(X_0, \dots , X_{2i-1}, \overset{2i}{0}, \overset{2i+1}{0}, 
                   X_{2i+1}, \dots , X_n) =0.
\end{equation}
We want to show $H(X_0, \dots, X_n)$
is the zero polynomial.
Since
the condition (\ref{conditionH}) holds for any $i$
with $1 \leq i \leq \frac{n-1}{2}$ and $X_0X_1$ divides $H$,
each monomial $X_0^{e_0}X_1^{e_1}\cdots X_n^{e_n}$
appeared in $H$ satisfies the condition
\begin{equation}\label{explicitconditionH}
\left\{
 \begin{split}
  & e_0 + \dots + e_n =\sqrt{q}+1 \\
  & e_0 > 0, \quad  e_1>0 \\
  & e_{2i} + e_{2i+1} >0 \quad \text{for $i$ with} \quad
     1 \leq i \leq \frac{n-1}{2}.
 \end{split}
\right.
\end{equation}
If $\sqrt{q}+1 < 2 + \frac{n-1}{2}$,
then no $(e_0, e_1, \dots , e_n)$
satisfies (\ref{explicitconditionH}).
Hence, in this case, $H$ is already the zero polynomial.

So we handle the opposit case below.
Put
\[
H(X_0, \dots , X_n) = \sum_{\mathbf{e}} c_{\mathbf{e}}
 X_0^{e_0}X_1^{e_1}\cdots X_n^{e_n},
\]
where $\mathbf{e} = (e_0, \dots , e_n)$
runs over the set of integer vectors
satisfying (\ref{explicitconditionH}).

Let
$\zeta$ be a root of $t^{\sqrt{q}-1}=-1$,
which is an element of $\mathbb{F}_q$.
Take a pair of nonsingular points in $X(\mathbb{F}_q)$
such a way that
\[
Q = (0, \dots , 0, \overset{2i}{1}, \overset{2i+1}{\zeta}, 0, \dots ,0)
\quad \text{and} \quad
Q' = (0, \dots , 0, \overset{2i}{\zeta}, \overset{2i+1}{1}, 0, \dots ,0).
\]
Since $\sqrt{q}+1 \geq 2 + \frac{n-1}{2} \geq 4$, 
$\sqrt{q}-1 \geq 2$. Also the characteristic of $\mathbb{F}_q$
and $\sqrt{q}-1$ are co-prime, we know $Q \neq Q'$.

Since $T_QX =\{-\zeta X_{2i} + X_{2i+1}=0 \}$
and $T_{Q'}X =\{ X_{2i} -\zeta X_{2i+1}=0 \}$,
we can apply Proposition~\ref{inductionstepZ} to this situation.
Especially, $X\cap T_QX$ is a cone of a Hermitian hypersurface.
Therefore 
\[
H(X_0, \dots , X_{2i}, \overset{2i+1}{\zeta X_{2i}},
X_{2i+2}, \dots , X_n)
\]
consists of Hermitian molecules.
Write down this polynomial explicitly:
\begin{multline*}
H(X_0, \dots , X_{2i}, \overset{2i+1}{\zeta X_{2i}},
X_{2i+2}, \dots , X_n)   \\
=\sum c_{\mathbf{e}} \zeta^{e_{2i+1}}X_0^{e_0}\cdots X_{2i-1}^{e_{2i-1}}
                          X_{2i}^{e_{2i}+e_{2i+1}}X_{2i+2}^{e_{2i+2}}
                            \cdots X_n^{e_n}  \\
= \sum_{\mathbf{e}'} \Bigl(
\sum_{v=0}^{\alpha} 
    c_{(e_0, \dots , e_{2i-1}, \alpha -v, v, e_{2i+2}, \dots , e_n)} \zeta^{v}
        \Bigr)X_0^{e_0}\cdots X_{2i-1}^{e_{2i-1}}
                          X_{2i}^{\alpha}X_{2i+2}^{e_{2i+2}}
                            \cdots X_n^{e_n},
\end{multline*}
where $\mathbf{e}'$ is the abbreviation for a $(n-1)$-pl
$(e_0, \dots , e_{2i-1}, e_{2i+2}, \dots , e_n)$
in $(e_0, \dots , e_{2i-1}, \alpha -v, v, e_{2i+2}, \dots , e_n)$.
Hence, for a fixed  
$\mathbf{e}'$,
\begin{equation}\label{coefficients0}
\sum_{v=0}^{\alpha} 
  c_{(e_0, \dots , e_{2i-1}, \alpha -v, v, e_{2i+2}, \dots , e_n)} \zeta^{v} =0
\end{equation}
for any $(\sqrt{q}-1)$-root $\zeta$ of $-1$.
Since $ \alpha \leq \sqrt{q}+1 -(2 + \frac{n-3}{2} ) < \sqrt{q}-1$,
all coefficients of $\zeta^{v}$ in (\ref{coefficients0})
are $0$.
Hence $H$ is the zero polynomial.
Therefore
\[
F = 
 X_0X_1 
  \bigl( \lambda X_0^{\sqrt{q}} + \lambda^{\sqrt{q}} X_1^{\sqrt{q}} \bigr)
  +
\sum_{i=1}^{\frac{n-1}{2}}
\Bigl(
X_{2i}^{\sqrt{q}}X_{2i+1} + X_{2i}X_{2i+1}^{\sqrt{q}},
\Bigr),
\]
which means $X$ is a Hermitian hypersurface.
Since $P_0 = (0,1,0, \dots , 0)$ is a nonsingular point of $X$,
$\lambda \neq 0$ by (\ref{pdF}).
Hence $X$ is nonsingular. 
\qed

\section{Supplementary}
In this section, we give two supplementaries.
\subsection{Comparison with Koen Thas' bound}
In \cite{tha2010}, Thas already gave another bound for $N_q(\mathcal{X})$
involving the invariant $k_{\mathcal{X}}$, where
$\mathcal{X}$ is a hypersurface of $\mathbb{P}^m$ of degree $d$
over $\mathbb{F}_q$ with $k_{\mathcal{X}} =k$.
Suppose $1 \leq k \leq m-2$. Then he proved that
\begin{equation}\label{KTbound}
N_q(\mathcal{X}) \leq dq^{m-1} + \theta_q(m-2)
              +(d-(q+1)) \sum_{i=k}^{m-2} q^i
                   \frac{\theta_q(m-1)}{\theta_q(i)\theta_q(i+1)}.
\end{equation}

\begin{proposition}
For $d$ with $d \leq q+1$,
the bound {\rm (\ref{GEB})} is better than {\rm (\ref{KTbound})}.
\end{proposition}
\proof
Let $S$ and $T$ be the upper bounds in (\ref{GEB}) and (\ref{KTbound})
respectively, namely,
\[
S= \theta_q(m-k-1)\cdot q^k(d-1) + \theta_q(k)
\]
and
\[
T=dq^{m-1} + \theta_q(m-2)
              +(d-(q+1)) \sum_{i=k}^{m-2} q^i
                   \frac{\theta_q(m-1)}{\theta_q(i)\theta_q(i+1)}.
\]
The claim is $T-S>0$ if $d \leq q+1$ and $1 \leq k \leq m-2$.
It is easy to see that
\[
S = \theta_q(m-1) +q^k + (d-2)q^k \theta_q(m-k-1)
\]
and
\[
T = \theta_q(m-1) + (d-1)q^{m-1} 
                   +(d-(q+1)) \sum_{i=k}^{m-2} q^i
                   \frac{\theta_q(m-1)}{\theta_q(i)\theta_q(i+1)}.
\]
Hence
\begin{equation}\label{TminusS1}
T-S = q^{m-1} -(d-2)q^k \theta_q(m-k-2)  - q^k
           +(d-(q+1)) \sum_{i=k}^{m-2} q^i
                   \frac{\theta_q(m-1)}{\theta_q(i)\theta_q(i+1)}.
\end{equation}
Let $t = q+1 -d$, which is nonnegative in the range of $d$.
Then the second term of the right-hand side of (\ref{TminusS1})
is rewritten as
\[
-q^{k+1} \theta_q(m-k-2) +(t+1)q^k \theta_q(m-k-2).
\]
Hence
\begin{multline}\label{TminusS2}
T-S = \\
q^{m-1} -q^{k+1} \theta_q(m-k-2) +(t+1)q^k \theta_q(m-k-2)
     -q^k  \\
      -t \sum_{i=k}^{m-2} q^i
                   \frac{\theta_q(m-1)}{\theta_q(i)\theta_q(i+1)}.
\end{multline}
Furthermore, since
\[
q^{m-1} -q^{k+1} \theta_q(m-k-2) =  -q^{k+1} \theta_q(m-k-3)
\]
and
\[
 -q^{k+1} \theta_q(m-k-3) + q^k \theta_q(m-k-2) -q^k =0,
\]
(\ref{TminusS2}) becomes
\[
T-S =
t \Bigl(
 q^k \theta_q(m-k-2) 
                   -\sum_{i=k}^{m-2} q^i
                   \frac{\theta_q(m-1)}{\theta_q(i)\theta_q(i+1)}
\Bigr).
\]
Since
\[
\frac{q^{i+1}}{\theta_q(i)\theta_q(i+1)} 
= \frac{1}{\theta_q(i)}-\frac{1}{\theta_q(i+1)},
\]
we get
\begin{align*}
T-S &= t \Bigl(
 q^k \theta_q(m-k-2) - \frac{\theta_q(m-1)}{q} 
     \bigl(\frac{1}{\theta_q(k)}-\frac{1}{\theta_q(m-1)} \bigr)
       \Bigr) \\
     &= \frac{t}{q\theta_q(k)}
     \bigl(
     q^{k+1} \theta_q(m-k-2)\theta_q(k) - \theta_q(m-1) + \theta_q(k)
     \bigr) \\
     & > \frac{t}{q\theta_q(k)}
        \bigl(
        q^{k+1} \theta_q(m-k-2) - \theta_q(m-1) + \theta_q(k)
        \bigr) =0.
\end{align*}
This completes the proof.
\qed

\subsection{The case where $m$ is even}
In Corollary~\ref{boundfornonsingular},
we gave an upper bound for $N_q(\mathcal{X})$
even if $\mathcal{X}$ is a nonsingular hypersurface
in an even dimensional projective space $\mathbb{P}^m$.
However, no nonsingular hypersurface achieves this upper bound
if $m$ is even.
More precisely, we can say:
\begin{annotation}
Suppose $m$ is even.
Let $\mathcal{X}$ be a hypersurface of degree $d \geq 2$ of
$\mathbb{P}^m$ over $\mathbb{F}_q$
with $ k_{\mathcal{X}} \leq \frac{m}{2}-1$. Then
\[
N_q(\mathcal{X}) \leq \theta_q(\frac{m}{2})q^{\frac{m}{2}-1}(d-1)
                           + \theta_q(\frac{m}{2}-1),
\]
however, equality no longer occurs.
\end{annotation}
\proof
This inequality comes from Theorem~\ref{theoremGEB}
like Corollary~\ref{boundfornonsingular} (ii) did.
Suppose equality holds for  $\mathcal{X}$.
Consider the ambient space $\mathbb{P}^m$ as
a hyperplane of $\mathbb{P}^{m+1}$,
and take $P_0 \in \mathbb{P}^{m+1} \setminus \mathbb{P}^m.$
Let $\tilde{\mathcal{X}} = P_0 \ast \mathcal{X} $
in $\mathbb{P}^{m+1}$.
Then
$\deg \, \tilde{\mathcal{X}} = \deg \, \mathcal{X},$
$k_{\tilde{\mathcal{X}}}  = k_{\mathcal{X}} +1$
and
\begin{align*}
N_q(\tilde{\mathcal{X}}) &= N_q(\mathcal{X})q +1 \\
             &= \theta_q(\frac{m}{2})q^{\frac{m}{2}}(d-1) +
                 \theta_q(\frac{m}{2}-1)q +1 \\
             &= \theta_q(\frac{m}{2}) \Bigl(
               (d-1) q^{\frac{m}{2}}+1 \Bigr).
\end{align*}
Let $n=m+1$. Then $\tilde{\mathcal{X}}$ satisfies the all assumptions 
of Theorem~\ref{actualtheorem} and equality holds in (\ref{bound}).
But from the latter part of this theorem, $\tilde{\mathcal{X}}$
must be nonsingular, which is a contradiction.
\qed

\gyokan

Finally we propose a conjecture for the case where $m$ is even.

\begin{conjecture}
Suppose $m \, ( \geq 4)$ is an even integer.
If $\mathcal{X}$ is a nonsingular hypersurface of degree $d$
in $\mathbb{P}^m$ over $\mathbb{F}_q$. Then
\[
N_q(\mathcal{X}) \leq \theta_q(\frac{m}{2}-1)
                      \Bigl(
               (d-1) q^{\frac{m}{2}}+1 \Bigr)
\]
might hold.
\end{conjecture}
When $m=2$, this inequality is just the Sziklai bound
and holds with only one exception \cite{hom-kim2010b}.
The nonsingular parabolic quadric hypersurface in $\mathbb{P}^m$,
and the nonsingular Hermitian hypersurface in $\mathbb{P}^m$
are examples for each of which equality holds.

\end{document}